\documentclass[11pt]{amsart}
\usepackage{eucal,amsmath,amssymb}
\usepackage[all]{xy}
\usepackage{hyperref}
\begin{document}

\parindent=0pt
\parskip=6pt

\newcommand{\ie}{\textit{i}.\textit{e}.\,}
\newcommand{\eg}{\textit{e}.\textit{g}.\,}
\newcommand{\cf}{{\textit{cf} \,}}
\newcommand{\que}{{\mathord{?}}}

\newcommand{\sslash}{{/\!\!/}}
\newcommand{\B}{{\mathbb B}}
\newcommand{\A}{{\mathcal A}}
\newcommand{\Ab}{{\mathbb A}}
\newcommand{\Bc}{{\mathcal B}}
\newcommand{\pB}{{\mathbb P}}
\newcommand{\C}{{\mathbb C}}
\newcommand{\Cc}{{\mathcal C}}
\newcommand{\Cp}{{\mathbb{C}_+}}
\newcommand{\D}{{\mathbb D}}
\newcommand{\bDel}{{\boldsymbol{\Delta}}}
\newcommand{\F}{{\mathbb F}}
\newcommand{\G}{{\mathcal G}}
\newcommand{\N}{{\mathbb N}}
\newcommand{\Oh}{{\mathbb O}}
\newcommand{\Q}{{\mathbb Q}}
\newcommand{\R}{{\mathbb R}}
\newcommand{\T}{{\mathbb T}}
\newcommand{\Z}{{\mathbb Z}}
\newcommand{\bu}{{\boldsymbol{\beta}}}
\newcommand{\fb}{{\mathfrak b}}
\newcommand{\x}{{\bf x}}
\newcommand{\z}{{\bf z}}
\newcommand{\bO}{{\bf{0}}}
\newcommand{\half}{{\textstyle{\frac{1}{2}}}}
\newcommand{\halfi}{{\textstyle{\frac{i}{2}}}}
\newcommand{\Con}{{\rm Con}}
\newcommand{\OCon}{{\rm OCon}}
\newcommand{\tCon}{{\widetilde{\rm Con}}}
\newcommand{\cC}{{\mathcal{C}}}
\newcommand{\tC}{{\widetilde{\mathcal{C}}}}
\newcommand{\Pic}{{\mathcal{P}ic}}
\newcommand{\fc}{{\mathfrak{c}}}
\newcommand{\dR}{{\rm dR}}
\newcommand{\tphi}{{\tilde{\phi}}}
\newcommand{\g}{{\mathfrak g}}
\newcommand{\w}{{\mathfrak w}}
\newcommand{\lines}{{\rm lines}}
\newcommand{\tot}{{\rm tot}}
\newcommand{\ab}{{\rm ab}}
\newcommand{\be}{{\bf e}}
\newcommand{\bk}{{\bf k}}
\newcommand{\bl}{{\bf l}}
\newcommand{\bv}{{\otimes}}
\newcommand{\mor}{{\rm Mor}}
\newcommand{\Gl}{{\rm Gl}}
\newcommand{\Sl}{{\rm Sl}}
\newcommand{\SF}{{\rm SF}}
\newcommand{\SU}{{\rm SU}}
\newcommand{\Hom}{{\rm Hom}}
\newcommand{\Fns}{{\rm Fns}}
\newcommand{\vt}{{\tau}}
\newcommand{\tr}{{\rm tr}}
\newcommand{\vT}{{\tilde{\tau}}}
\newcommand{\ve}{{\varepsilon}}
\newcommand{\Arg}{{\rm Arg \:}}
\newcommand{\Pc}{{\mathcal P}}
\newcommand{\Pm}{{\mathbb P}}
\newcommand{\Rp}{{\rm Re \:}}
\newcommand{\bD}{{\boldsymbol{\Delta}}}
\newcommand{\Arctan}{{\rm Arctan \:}}
\newcommand{\sgn}{{\rm sgn}}
\newcommand{\bpi}{{\boldsymbol{\pi}}}
\newcommand{\Sc}{{\mathcal S}}
\newcommand{\String}{{\rm String}}
\newcommand{\Hs}{{\sf H}}
\newcommand{\Mod}{{\rm Mod}}
\newcommand{\Maps}{{\rm Maps}}
\newcommand{\pic}{{\rm Pic}}

\newcommand{\SB}{{{\rm S}\mathbb{B}}}
\newcommand{\tN}{{\widetilde{\mathbb{N}}}}
\newcommand{\ph}{{\hat{p}}}
\newcommand{\Spin}{{\rm Spin}}
\newcommand{\alg}{{\rm alg}}
\newcommand{\bs}{{\boldsymbol{\sigma}}}
\newcommand{\bK}{{\bf K}}
\newcommand{\ord}{{\rm ord}}
\newcommand{\op}{{\rm op}}
\newcommand{\hp}{{\hat{p}}}
\newcommand{\Top}{{\rm top}}
\newcommand{\tbeta}{{\tilde{\beta}}}
\newcommand{\Tof}{{{\mathcal T}_{\rm OF}}}
\newcommand{\dvol}{{d{\rm vol}}}
\newcommand{\tLG}{{\widetilde{\rm{L}G}}}
\newcommand{\SO}{{\rm SO}}
\newcommand{\bplus}{{\boldsymbol{+}}}
\newcommand{\tfb}{{\tilde{\mathfrak b}}}
\newcommand{\syl}{{\sf s}}
\newcommand{\ess}{{\sf S}}
\newcommand{\pp}{{\: +_T \:}}
\newcommand{\tp}{{\: \tilde{+}_T \:}}
\newcommand{\e}{{\bf e}}
\newcommand{\wz}{{\bf wzw}}
\newcommand{\wzw}{{\mathfrak{wzw}}}
\newcommand{\gl}{{\mathfrak{gl}_1}}
\newcommand{\Tor}{{\rm Tor}}
\newcommand{\Or}{{\rm or}}

\title{Toward the group completion of the Burau representation}

\author[Jack Morava]{Jack Morava}

\address{Department of Mathematics, The Johns Hopkins University,
Baltimore, Maryland 21218}

\email{jack@math.jhu.edu}

\author[Dale Rolfsen]{Dale Rolfsen}

\address{Pacific Institute for the Mathematical Sciences and Department of 
Mathematics, University of British Columbia, Vancouver, BC, Canada V6T 1Z2}

\email{rolfsen@math.ubc.ca}

\begin{abstract} {Following Boardman-Vogt, McDuff, Segal, and others, we construct a monoidal topological groupoid or space of finite subsets of the plane, and interpret the Burau representation of knot theory as a topological quantum field theory defined on it. Its determinant or {\bf writhe} is an invertible braided monoidal TQFT which group completes to define a  Hopkins-Mahowald model for integral homology as an $E_2$ Thom spectrum. We use these ideas to construct an infinite cyclic (Alexander) cover for the space of finite subsets of $\C$, and we argue that the TQFT defined by Burau is closely related to the SU(2)-valued Wess-Zumino-Witten model for string theory on $\R^3_+$.}\end{abstract}

\maketitle \bigskip

{\bf \S 0} Introduction

{\bf  \S 1} Topological categories and their geometric realizations

 {\bf \S 2} Monoidal categories of finite subsets of the plane

{\bf \S 3} Categorification of the Burau representation

{\bf \S 4} The braid monoid $\B_*$

{\bf \S 5} The restricted braid monoid $S\B_*$

{\bf \S 6} Some calculations in Artin coordinates

{\bf \S 7} Monoidal braidings and $p$-adic spheres

{\bf \S 8} The Hopkins-Mahowald theorem $H\Z \simeq M\SF(3)$  \bigskip

{\bf \S 0} {\sc Introduction}

The classical {\bf un}reduced Burau representations
\[
\B_k \to \Gl_k(\Z[t,t^{-1}])
\]
of Artin's braid groups send the juxtaposition or $\otimes$ composition of braids to direct sum $\oplus$ of matrices, defining a functor 
\[
\fb_* : [* \sslash \B_*] \to [* \sslash \Gl_*(\Lambda)] \;,
\]
between monoidal categories (with $\Lambda := \Z[t,t^{-1}]$); the geometric realization or classifying space construction \cite{54} then defines a morphism 
\[
\coprod_{k \geq 0}  B\B_k \to \coprod_{k \geq 0}B\Gl_k(\Lambda)
\]
of topological monoids. The group completion morphism $M \to M^+ := \Omega BM$ associates to this, a morphism
\[
\Omega^2 S^2 \to \Z \times \Omega^\infty K(\Lambda)
\]
of ($E_2$) loop-spaces, the latter term representing Quillen's algebraic $K$-theory of $\Lambda$. The determinant construction $K(\Lambda) \to \pic(\Lambda)$ (which sends a finitely generated free $\Lambda$-module to its top exterior power regarded as a $\Lambda$-line), together with a choice $t \mapsto -a \in \Z_p^\times$ and the Hopf fibration $S^3 \to S^2$, define a morphism
\[
\Z \times \Omega^2 S^3_p \to \pic (\Z_p)
\]
of ($p$-complete $E_2$) loop spaces. Recent work \cite{5, 18} of Clausen and of Beaudry, Goerss, Hopkins, and Stojanoska extends this construction to define a morphism
\[
 |[* \sslash \B_*]|^+_p \simeq \Z \times \Omega^2 S^3_p \to \Pic_@(S^0_p) \simeq \Z \times B\Sl(S^0_p)
  \]
of braided monoids, with the latter term referring to the Picard group of oriented $p$-adic spheres under smash product. The map induced on universal covers defines a family of $S^0_p$-lines, and thus \cite{2} a generalized Thom spectrum which we propose to call $M\SF(3)$. In \S 8 we argue that this spectrum is a model for Hopkins and Mahowald's \cite{42, 43} construction of the integral Eilenberg-Mac Lane spectrum $H\Z$.

The isomorphism $\B_k \cong \pi_0\D_k$ of the braid groups with mapping class groups of punctured spheres leads to an interpretation of the determinant of the Burau representation as a topological quantum field theory for braids, regarded as particle trajectories in 3-space. We propose to interpret the restriction of this construction to the group completion of the submonoid of braids with classical writhe zero as a homotopy-theoretic TQFT defined on particles with trajectories subject to moves of Reidemeister type II and III but not type I, \ie without kinks.  \newpage

{\bf  \S 1} {\sc Topological categories and their geometric realizations} 

We work in a semiclassical language of topological categories \cite{6,33,58} with compactly generated spaces of objects and morphisms\begin{footnote}{We write $=$ for identity, $:=$ for definitions, $\cong$ for isomorphism, and $\simeq$ for homotopy equivalence of spaces or equivalence of categories. Groups act on the left but we write their quotients on the right.}\end{footnote}. The transformation groupoid $[X \sslash G]$ associated to a continuous action 
\[
\alpha : G \times X \to X
\]
of a topological monoid, \eg a group $G$, on a space $X$ is the motivating example: the elements $x_0,x_1 \in X$ are its objects, and the `transporter'
\[
\{g \in G \:|\: \alpha(g,x_0) = x_1 \} = \mor_{X \sslash G}(x_0,x_1)
\]
is the space of morphisms from $x_0$ to $x_1$. A space $X$ can be regarded as a topological category $[X \sslash 1]$ with only identity morphisms, and a topological group (or monoid) $G$ can be regarded as the category $[* \sslash G]$ with only one object. A morphism $(G,X) \to (H,Y)$ of group actions becomes a functor $[X \sslash G] \to [Y \sslash H]$ of categories. 

Following Grothendieck and Segal \cite{54} a topological category $\A$, with objects $\A[0]$ and morphisms $\A[1]$ defines a simplicial space 
\[
* \leftarrow \A[0] \Leftarrow \A[1] \Lleftarrow A[2] = \A[1] \times_{\A[0]} \A[1] \cdots
\]
(with $\A[k]$ as the space of $k$-tuples of composable morphisms); there is then a product-preserving `geometric realization' functor
\[
\A \mapsto |\A| := (\coprod_{k \geq 0} \A[k] \times \Delta^k)/({\rm relations)}
\]
to topological spaces. This totalization provides a canonical pointed classifying space $BG := |[* \sslash G]|$ for a topological monoid: $B\{1\}$, for example, is a kind of thick point. We have a fibration
\[
X \to |[X \sslash G]| \simeq EG \times_G X \to EG \times_G * \simeq BG
\]
which identifies $|[X \sslash G]| \simeq X_{hG}$ with the Borel construction or homotopy quotient of $X$ by $G$; the functor 
\[
[X \sslash G] \to [(X/G) \sslash 1]
\]
is then a kind of resolution of the topological quotient. More generally, the group completion theorem \cite{46} asserts that the corresponding construction for a monoid is a homology fibration.

Note for example that 
\[
|[X \sslash 1]| \simeq B\{1\} \times X \to X
\]
is a homotopy equivalence but not generally a homeomorphism. A closed subgroup $H \subset G$ defines an equivalence $[G \sslash H] \to [(G/H) \sslash 1]$ of categories, and the equivalence 
\[
\xymatrix{
\G \ar[r]^-\simeq & \coprod_{c \in \pi_0\G} {\rm Aut_\G}(c)}
\]
of a discrete groupoid $\G$ with its skeleton extends to imply a homotopy equivalence of geometric realizations
\[
|[(G/H) \sslash G]| \simeq |[* \sslash H]| \;,
\]
in the topological case. \bigskip

 {\bf \S 2} {\sc Monoidal categories of finite subsets of the plane} \bigskip

{\bf 2.1} Disjoint union defines a monoidal structure on the category of  finite sets, but for subsets of a given space things are more complicated. In the theory of configuration spaces, the symmetric group $\Sigma_k$ acts on a $k$-fold product $Z^k$ by permuting coordinates, preserving the `thick diagonal'
\[
\bD^k(Z) = \{\z = (z_1,\dots,z_k) \in Z^k \:|\: \exists i \neq j \; {\rm such \;; that} \; z_i = z_j \}
\]
consisting of sequences with at least one repetition. Thus $\bD^k(\C)$ is the hyperplane arrangement defined by the union of the kernels of the linear maps
\[
\C^k \ni \z \mapsto z_i - z_j \in \C
\]
($i \neq j$), and $\Sigma_k$ acts freely on its complement, the space $\OCon^k(\C) := \C^k -\bD^k(\C)$ of ordered (\ie labelled) $k$-tuples of distinct points, with a quotient
\[
\OCon^k(\C)/\Sigma_k := \Con^k(\C) 
\]
which can be regarded as the space of subsets of cardinality $k$ of the plane. We will write $\{z\}$ or even $z$ for the equivalence class or configuration defined by $\z \in \C^k$; in particular we will write $\{k\}$ for the set $\{1,\dots,k\}$ defined by a natural number $k$, and take it to be the basepoint $*_k$ for $\Con^k(\C)$. 

We thank Alan Hatcher for reminding us that classical work of Smale \cite{57} implies that the topological group $\D$ of compactly supported diffeomorphisms of the plane is contractible. For our purposes it will be useful to interpret $\D$ as the group of diffeomorphisms of the Riemann sphere 
\[
\R^2_+ = \C_+ = \C P^1 \cong S^2 
\]
which leave fixed a neighborhood of the point at infinity. The diagonal action defines a transitive group action
\[
\D \times \Con^k(\C) \to \Con^k(\C)
\]
implying a diffeomorphism 
\[
\xymatrix{
\D/\D_{\{k\}} \ar[r]^\cong & \Con^k(\C) \;,}
\]
where $\D_{\{k\}}$ is the isotropy group or stabilizer of the configuration $\{k\}$. We can thus regard 
\[
\Cc^k := [(\D/\D_{\{k\}}) \sslash \D] \cong [\Con^k(\C) \sslash \D]
\]
as the topological groupoid of subsets of the plane with cardinality $k$, with compactly supported diffeomorphisms as morphisms, and we will write
\[
\Cc^* := \coprod_{k \geq 0} \Cc^k
\]
(with $\Cc^0$ a singleton) for the corresponding category of finite subsets of $\C$. The space
\[
\Cc(\{w\},\{z\}) := \{ \phi \in \D \:|\: \phi \{w\} = \{z\} \}
\]
of morphisms in $\Cc_*$ is empty unless $\{z\}$ and $\{w\}$ have the same cardinality. 

{\bf Proposition 1} :  {\it Geometric realization of the categories and functors in the composition
\[
[\Con(\C)^* \sslash 1] \cong [(\D/\D_*) \sslash 1] \to [(\D/\D_*) \sslash \D] = \Cc^* \to [* \sslash \D_*]
\]
defines a sequence of homotopy equivalences.} 

The final arrow, for example, sends a morphism $\phi : \D_k \gamma_0^{-1} \to \D_k \gamma_1^{-1}$ to $\gamma_1^{-1} \phi \gamma_0 \in \D_k$. These categories thus differ, but are nevertheless in some sense equivalent geometrically. Note that this composition of morphisms is directed, rather than being a zigzag equivalence.  \bigskip

 {\bf Remark} We take as known the profound fact that the diffeomorphism groups $\D_{\{k\}}$ (and hence the groups $\D_{[\bk]}$ defined below) are homotopically discrete, making the homomorphisms
\[
\D_{\{k\}} \to \pi_0 \D_{\{k\}} \cong \B_k
\]
to the mapping class groups of the punctured plane into homotopy equivalences \cite{7}(\S 1.3), and that these groups of components can be identified with Artin's braid groups $\B_k$. We will return to this in \S 4. \bigskip

{\bf 2.2}  Let us recall Segal's variant \cite{55} of Boardman and Vogt's \cite{10} construction of a rigidified monoidal category $\Cc^\otimes \to \Cc^*$ of unordered points, together with a forgetful functor to $\Cc^*$ which defines an equivalence on geometric realizations.
 
{\bf Definition} If $k \in \N_{\geq 1}$ let
\[
[\bk] := \{ z \in \C \:|\: \half \leq {\rm Re} \: z \leq k + \half \} \cong [\half,k+\half] \times i\R
\]
(but $[{\bf 0}] = \emptyset$), and let $\Con^{\bk}(\C) \subset \Con^k(\C)$ denote the space of cardinality $k$ subsets of the {\bf interior} of $[\bk]$; the inclusion is a homotopy equivalence. If $k' \in \N_{\geq 0}$, let
\[
z \mapsto T_{k'}(z) = z + k' :[\bk] \to T_{k'}[\bk] \subset [\bk + \bk']
\]
be the translation of $[\bk]$ by $k'$ steps to the right; this then defines a commutative monoid $\{[\bk] \:|\: k \in \N_{\geq 0} \}$ of subsets of $\C$, with unit $[\bO] := \emptyset$ and associative composition
\[
[\bk] + [\bk'] = [\bk] \cup T_{k}[\bk'] \;:
\]
\[
([\bk] + [\bk']) + [\bk''] = [\bk + \bk'] + [\bk''] = [\bk + \bk' + \bk''] = [\bk] + [\bk' + \bk''] = [\bk] + ([\bk'] + [\bk'']) \;.
\]
This is compatible with the ordinal sum composition on the collection $\{ \{k\} \:|\: k \in \N \}$ of basepoints for our configuration spaces, defined by regarding $k$ as $\{k\} \times 0 \subset [\bk]$. Mac Lane's pentagon condition \cite{34} for a monoidal category holds for this composition: it asserts that composing these operations in two different ways defines the same map
\[
(([\bk] + [\bk']) + [\bk'']) + [\bl] \to [\bk] +(([\bk'] + [\bk'']) + [\bl])
\; .
%  \; (= [\bk + \bk' + \bk'' + \bl])
\] 

Let $\D_{[\bk]} \subset \D$ be the subgroup of diffeomorphisms supported in the interior of $[\bk]$; similarly let $\D_{\{\bk\}} := \D_{[\bk]} \cap \D_{\{k\}}$ be the group of diffeomorphisms of the `time zone' $[\bk]$ supported in its interior and fixing $\{k\}$. We then have a commutative diagram 
\[
\xymatrix{
\D_{\{\bk\}} \ar[d]^\simeq \ar[r]^-\subset & \D_{[\bk]} \ar[d]^\simeq \\
\D_{\{k\}} \ar[r]^-\subset & \D }
\]
of topological categories and functors
\[
\Cc^{\otimes k} := [(\D_{[\bk]}/\D_{\{\bk\}}) \sslash \D_{[\bk]}] \to [(\D/\D_{\{k\}}) \sslash \D] = \Cc^k
\]
making $\Cc^\otimes := \coprod_{k \geq 0} \Cc^{\otimes k}$ into a monoidal category with geometric realization $|\Cc|^*$, with a (not so commutative) composition
\[
\pp : \Con^{\bk'}(\C) \times \Con^{\bk''}(\C) \to \Con^{\bk' + \bk''}(\C)
\]
defined by
\[
\{w\} \;,\:\{z\} \mapsto \{w\} \cup T_{k'}\{z\} = \{w_i, z_j + k'\} := \{w \pp z\} \;. 
\]

{\bf 2.3} The discriminant
\[
\Delta\{z\} = \prod_{i \neq j} (z_i - z_j)  = (-1)^{\binom{k}{2}} \prod_{i < j} (z_i - z_j)^2 \in \C^\times \cong \R^\times_+ \times \T
\]
of $\{z\} \subset \C$ defines smooth translation-invariant functions
\[
\Delta : \Con^k(\C) \to \C^\times \,
\]
\eg $\Delta (*_k) =  \Delta \{k\} =  (-1)^{\binom{k}{2}} (\prod_{1 \leq j \leq k} \Gamma(k))^2$ is a kind of superfactorial. Let $\Theta = \Arg \Delta \in \T \cong \R/\Z$ denote its phase or complex argument.

In 1969 F Gorin and V Lin \cite{32} showed that 
\[
\Delta_* : \Cc^* \to \N \times \C^\times
\] 
is a locally trivial fibration, and therefore much more; see \S 5.2, 6.3 below. Their result implies that the pullback 
\[
\xymatrix{
\tCon^k(\C) \ar@{.>}[d] \ar@{.>}[r] & \C \ar[d]^\e \\
\Con^k(\C) \ar[r]^-\Delta & \C^\times }
\]
(with $\e(t) = \exp(2 \pi i t)$) defines an infinite cyclic cover of the configuration space $\Con^*(\C)$. Its elements are pairs $(\{u\},\theta)$ with $\theta \in \R$ such that $\Theta \{u\} \equiv \theta$ modulo $\Z$, with elements $\phi \in \D$ acting by 
\[
\phi, \: \{z,\alpha\} \mapsto \{\phi\{z\},\; \alpha + \Theta(\phi\{z\}) - \Theta(\{z\})\} \;.
\]

{\bf Proposition 2} :  {\it Disjoint union $\pp$ on $\Cc^\otimes$ lifts
\[
\xymatrix{
\tCon^{\bk'}(\C) \times \tCon^\bk(\C) \ar[d] \ar[r]^-\tp & \tCon^{\bk + \bk'}(\C) \ar[d] \\
\Con^{\bk'}(\C) \times \Con^\bk(\C) \ar[r] \ar[r]^-\pp & \Con^{\bk + \bk'}(\C)}
\]
to an associative composition law 
\[
\{u,\theta\} \tp \{v,\tau\} = \{u \pp v, \theta + 2 \Arg \syl(u,v) + \tau \}
\]
on $\tC^\otimes$, where $\syl(u,v) = {\rm Res}(\{u\},T_{k'}\{v\}) := \prod (u_i - (v_j n + k'))$ is Sylvester's resultant.}

{\bf Proof:} Checking the commutativity is straightforward from the definitions. Perhaps more conceptually, the cocycle
\[
\syl(u,v)^2 = \frac{\Delta\{u \pp v\}}{\Delta\{u\} \cdot \Delta \{v\} } \in Z^2(\Con^\otimes(\C),\T) 
\]
is a coboundary, making the associativity lifting obstruction equal to zero in $H^1(\Cc^\otimes,\T)$. 

{\bf Corollary}
\[ 
\emptyset \times \Z  \rightarrowtail \tC^\otimes \; \twoheadrightarrow \; \Cc^\otimes
\]
{\it presents $\tC^\otimes$ as a central extension or abelian cover of the monoidal category} $\Cc^\otimes$. 

We propose to think of the objects of $\tC^\otimes$ as collections $(\{u\},\theta)$ of particles in the plane endowed with an intrinsic collective real-valued topological charge. \bigskip

{\bf \S 3} {\sc a categorification of the Burau representation} 

In this language the geometric construction of W. Burau \cite{7,8,14} defines a monoidal `quantum field theory' for finite subsets of $\C$ under compactly supported diffeomorphisms, with values in CW-spaces, taking disjoint union to wedge sum:

{\bf 3.1} If $\{z\} \in \Con^k(\C)$, the one-point compactification of $\C - \{z\}$ is homeomorphic to the space obtained by collapsing the subset $\{z+\} :=\{z\} \cup \{+\} \subset\C_+$ to a point; we will simplify notation by writing $\C\{z\} = \C_+/\{z+\}$ for this quotient space. Its fundamental group $\pi_1(\C\{z\},+) \cong F_k$ is free, generated by loops $\gamma_i : (\R_+,+) \to (\C\{z\},+)$ defined by an arc from $\infty$ to $z_i$ (a loop since $\infty$ and $z_i$ are identified); thus $H_1(\C\{z\},\Z)$ is free abelian of rank $k$, generated by the Hurewicz images $[\gamma_i]$, while
\[
H_2(\C\{z\},\Z) \cong H_2(\C_+,\{z+\};\Z)) \cong \Z \;.
\]
If $\tphi : \C\{z\} \to \C\{z'\}$ is the homeomorphism of one-point compactifications associated to $\phi \in \Cc(\{z\},\{z'\})$, then the Hurewicz homomorphism defines a commutative diagram
\[
\xymatrix{
\pi_1(\C \{z\},+) \ar[d] \ar[r]^{\pi(\tphi)} & \pi_1(\C \{z'\},+) \ar[d] \\
H_1(\C \{z\};\Z) \ar[d]^\tr \ar[r]^{H_1(\tphi)} & H_1(\C \{z\};\Z)
\ar[d]^{\tr'}\\
\Z \ar[r]^= & \Z \;.}
\]
The trace homomorphism in the bottom square sends $\sum a_i[\gamma_i]$ to $\sum a_i \in \Z$: $\D_k$  acts on the homology of $\C\{z\}$ through the quotient $\B_k \to \Sigma_k \to \Gl_k(\Z)$, defining homomorphisms
\[
H_1(\C\{k\}) \to H_0(\Sigma_k,H_1\C\{z\}) \cong \Z 
\]
and thus lifts 
\[
\xymatrix{
(B\{z\},*) \ar[d]^c \ar[r]^{B(\phi)} & (B\{z\},*) \ar[d]^{c'}\\
(\C\{z\},+) \ar[r]^\tphi & (\C\{z'\},+) }
\]
of $\phi$ to certain infinite cyclic ($\ie$ Burau) covering spaces: Riemann surfaces equivariant with respect to canonical actions of the deck-transformation group $\Z$. For example $\C\{1\} \cong \C P^1/\{1,\infty\}$, while $B\{1\}$ is essentially the Riemann surface of $\log z$. 

The space $V\{z\} := B\{z\}/c^{-1}(+)$ underlies the {\bf un}reduced Burau representation.  Since the homology groups in this section all have coefficients in $\Z$, from now on we will omit it from the notation. Let us denote the deck-transformation group ring $\Z[\Z]$ by $\Lambda = \Z[t,t^{-1}]$. Note that if $\phi \in \Cc(\{z\},\{z'\})$ and $\phi' \in \Cc(\{z'\},\{z''\})$, then
\[
B(\phi') \circ B(\phi) = B(\phi' \circ \phi)
\]
as maps equivariant with respect to the deck-transformation action.

{\bf 3.2}  The folding map $s_k : \Cp \vee \Cp \to \Cp$ (which sends $\{z\}$ in the first $\Cp$ of the domain to itself, and $\{z'\}$ in the second $\Cp$ to $\{k + z'\}$) defines a commutative diagram
\[
\xymatrix{
\Cp \vee \Cp \ar[d] \ar[r]^{s_k} & \Cp \ar[d] \\
\C\{z\} \vee \C\{z'\} \ar[r]^{\tilde{s}_k} & \C\{z+z'\} }
\]
of base-pointed spaces (with the vertical maps defined as collapses, as above).
If
\[
V\{z\} := B\{z\}/c^{-1}(+)
\]
then the commutative diagram
\[
\xymatrix{
\pi_1(\C\{z\} \vee \C\{z'\},+) \cong F_k * F_{k'} \ar[d] \ar[r]^-{\pi(s_k)} & \pi_1(\C\{k+k'\},+) \cong F_{k+k'} \ar[d] \\
H_1(\C\{z\} \vee \C\{z'\}) \cong \Z^k \oplus \Z^{k'}  \ar[d]^\tr \ar[r]^-{H(s_k)} & H_1(\C\{z+z'\}) \cong\Z^{k+k'} \ar[d]^{\tr'} \\
\Z \ar[r]^= & \Z }
\]
defines a lift to a map of covering spaces, from $B\{k\} \cup_{c^{-1}(+)} B\{z'\}$ over $\C\{z\} \vee \C\{z'\}$ to $B\{z+z'\}$ over $\C\{z+z'\}$, as well as maps
\[
V\{z\} \vee V\{z'\} \to V\{z+z'\} \;,
\]
equivariant with respect to the deck-transformation action. We have short exact sequences
\[
0 \to H_1(B\{z\}) \to H_1(B\{z\},c^{-1}(+)) \cong H_1(V\{z\} \to H_0(c^{-1}(+)) \cong \Lambda \to 0
\]
by Eilenberg - Steenrod, and isomorphisms
\[
0 = H_1(B\{z\} \cap B\{z'\},c^{-1}(+)) \to H_1(B\{z\},c^{-1}(+)) \oplus H_1(B\{z'\},c^{-1}(+)) =
\]
\[
H_1(V\{z\}) \oplus H_1(V\{z'\}) \cong H_1(V\{z+z'\}) =  H_1(B\{z+z'\},c^{-1}(+)) \to 0
\]
of (unreduced) Burau modules, free of rank $k+k'$ over $\Lambda$, by Mayer-Vietoris. We write $\tfb$ for the reduced Burau representation, defined by $H_1(B\{z\})$.  \bigskip

{\bf Proposition 3} : {\it The unreduced Burau representation
\[
\fb : \{z\} \mapsto \fb\{z\} := H_1(V\{z\},\Z) : (\Cc^\otimes, +_T) \to (\Lambda-\Mod,\oplus)
\]
is a product-preserving functor from the monoidal category $\Cc^\otimes$ to the groupoid of finitely generated free $\Lambda$-modules and their isomorphisms.} $\Box$ 

We return to these groups in \S 6. \bigskip

{\bf \S 4} {\sc The braid monoid} \bigskip

 {\bf 4.1} Following Joyal and Street, the collection
\[
\langle \sigma_i, \; 1 \leq i \leq k-1 \:|\: |i-k| > 1 \Rightarrow [\sigma_i,\sigma_k] = 1, \; \sigma_i \sigma_{i+1} \sigma_i = \sigma_{i+1} \sigma_i \sigma_{i+1} \rangle
\]
of braid groups $\B_k$ on $k > 0$ strings defines a free monoidal category \cite{34}, \cite{37} (XIII.3)
\[
\Bc_* := \coprod_{k \geq 0} [k \sslash\B_k]
\]
on one generator, with strictly associative tensor product
\[
\otimes: \B_k \times \B_{k'} \to \B_{k+k'}
\]
defined by juxtaposition of isotopy classes of braid diagrams. We take $\B_0$ and $\B_1$ to be singleton categories, each with one object and one morphism. The generator $\sigma_i$ is the class of a Dehn half-twist along the interval $[i,i+1]$, and the full $k$-strand twist $z_k = (\sigma_1 \dots \sigma_{k-1})^k$ generates the center of $\B_k$, with $\w(z_k) = k(k-1)$.
 
{\bf 4.2} Dehornoy \cite{23} discovered that braid groups are left-orderable, in the sense that there is a strict total ordering $<$ of its elements with the property that $h < h'$ implies $gh < gh'$ for elements $g, h, h'$. The ordering $<$ is defined as follows:  A braid is called $\sigma$-positive if it has an expression in the Artin generators such that the $\sigma_i$ with lowest index $i$ occurs with only positive exponent. The key to the ordering is the following nontrivial fact \cite{31}.

{\bf Proposition 4} :  {\it For a braid $\beta \in \B_k$ exactly one of the following holds: (1) $\beta$ is $\sigma$-positive, (2) $\beta^{-1}$ is $\sigma$-positive, (3) $\beta$ is the identity element $1_k \in \B_k$.}

{\bf Definition:} For braids $\alpha, \beta \in \B_k$ declare $\alpha < \beta$ iff $\alpha^{-1}\beta$ is $\sigma-$positive.

It is easily checked using the Proposition that this defines a left-ordering of $\B_k$.  This ordering is compatible with the canonical inclusions $\B_k \subset \B_{k+1}$.  We remark that it is also a left-ordering of the Garside braid monoid $\B_k^+$ consisting of all braids expressible as words in the Artin generators with only positive exponents, and is in fact a well-ordering there.

By \cite{30} the groups $\pi_0\D_n$ inherit Dehornoy's order: a mapping class of the complex plane fixing $\{n\}$ is considered positive if the image of the real axis under a suitable representative first departs from the real axis (going from left to right) into the lower half plane; this is the opposite convention of that of \cite{24}.

If $\alpha \in \B_k$ and $\beta \in \B_{k'}$, their tensor product $\alpha \otimes \beta$ may be expressed in terms of the Artin generators by taking the product of $\alpha$ with the word expressing $\beta$ with $k$ added to all the indices of the generators occurring in $\beta$. Since the resulting words are in subgroups $\langle \sigma_1, \dots , \sigma_{nk1}\rangle$ and $\langle \sigma_{k+1}, \dots , \sigma_{k + k'-1}\rangle$ of $\B_{k+k'}$,  respectively, they commute.  If $\alpha, \alpha' \in \B_n$ and $\beta, \beta' \in \B_{k'}$ we have $(\alpha \otimes \beta)(\alpha' \otimes \beta') = \alpha \alpha' \otimes \beta\beta'$ and $(\alpha \otimes \beta)^{-1} = \alpha^{-1} \otimes \beta^{-1}$. Moreover, if $\alpha$ is $\sigma$-positive, then so is $\alpha \otimes \beta$ for arbitrary $\beta \in \B_{k'}$.

{\bf Proposition 5} : {\it The lexicographic order on $\B_k \otimes \B_{k'}$ coincides with its Dehornoy order, regarded as a subset of $\B_{k+k'}$. That is, $\alpha \otimes \beta < \alpha' \otimes \beta'$ if and only if $\alpha < \alpha'$ or $\alpha = \alpha'$ and $\beta < \beta'$. Tensor product of braids is thus a monotone function.}

{\bf Proof:}  $\alpha \otimes \beta < \alpha' \otimes \beta' \iff (\alpha \otimes \beta)^{-1}(\alpha' \otimes \beta') = \alpha^{-1}\alpha' \otimes \beta^{-1}\beta'$ is 
$\sigma$-positive $\iff \alpha^{-1}\alpha'$ is $\sigma$-positive or $\alpha^{-1}\alpha' = 1_k$ and $\beta^{-1}\beta'$ is $\sigma$-positive. $\qed$

The twist class is positive in the Dehornoy ordering: the generator $\sigma_{k'+k-1}$ has the highest subscript occurring and it occurs only once (with positive exponent) so $\tau_{k,k'}$ is $\sigma_{k'+k-1}$ - positive and therefore greater than the identity. \bigskip

{\bf \S 5} {\sc The restricted braid monoid $S\B_*$}

{\bf 5.1} The braid relations imply that
\[
(\sigma_i \sigma_{i+1})\sigma_i(\sigma_i \sigma_{i+1})^{-1} = \sigma_{i+1}\sigma_i\sigma_{i+1}\sigma_{i+1}^{-1}\sigma_i^{-1} = \sigma_{i+1}
\]
and hence that the Artin generators are all conjugate to one another; it follows that they map to the same element of the abelianization
\[
\B_k/[\B_k,\B_k] := \B_k^\ab \cong \Z \;,
\]
defining homomorphisms
\[
\w_k : \B_k \to \B_k^\ab \cong \Z
\]
which send a word $b \in \B_k$ in the $\sigma_i$ to the sum $\w(b)$ of its Artin exponents. The braid closure $b^*$ of $b$ is a framed oriented link, whose {\bf writhe}, \ie its (invariant) number of overcrossings minus undercrossings, equals $\w(b)$. The commutative diagram
\[
\xymatrix{
\B_k \ar[d] \ar[r]^-\w & \Z \ar[d] \\
\Sigma_k \ar[r]^-\sgn & \Z_2 }
\]
suggests interpreting the writhe of a braid as a generalization of the sign of a permutation. Like the sign, it is additive under juxtaposition. 

{\bf 5.2} We are indebted to Juan Gonzales-Meneses for pointing us to work of Gorin and Lin \cite{32}(Lemma 3.6), which is perhaps not as well known as it could be:

{\bf Theorem} : {\it The homomorphism
\[
\Delta_{k*} : \pi_1(\Con_k(\C),\{k\}) \cong \B_k \to \pi_1(\C^\times,*_k) 
\cong \Z                                                                      
\]                                                                            
(induced on fundamental groups by the discriminant) equals the writhe $\w_k$.}

These authors show that $\Delta_k$ is a locally trivial fibre map whose fibre is connected and has as fundamental group the derived or commutator subgroup
\[
\SB_k = [\B_k,\B_k] = \ker \w_k \;;
\]
the homotopy exact sequence of the fibration then implies that $\Delta_{k*}$ equals the writhe up to sign. They show moreover that $\SB_k$ is perfect for $k>4$; if $k = 1,2$ these groups are trivial, while $\SB_3 \cong F_2$ is free of rank two (the trefoil complement has a punctured torus as fiber), and $\SB_4$ is a semidirect product of two copies of $F_2$. We take $\SB_0$ to be the empty monoid. By now, the homology of the fiber \cite{15,17,22,25,39} has an extensive literature. 

We present a proof by direct calculation below, which gives explicit conventions for signs. The discriminant
induces a homomorphism 
\[
\Delta_* : \B_k = \pi_1(\Con^k(\C),\{k\}) \to \pi_1(\C^\times,1) = \Z.
\]
\noindent
The {\em writhe} $\w(\beta)$ of a braid $\beta = \sigma_{k_1}^{e_1} \cdots \sigma_{k_m}^{e_m}$ is just the exponent sum $e_1 + \cdots + e_m$ when written in terms of the Artin generators.  We show the following. \bigskip

{\bf Proposition 6} : {\it For $\beta \in \B_k, \; \Delta_*(\beta) = \w(\beta)$}
 \bigskip

Proof : With $\{1, 2, \dots, k\}$ as the basepoint of $\Con^k(\C)$, the braid $\sigma_1$ may be represented by a $180^o$ rotation of the points $\{1, 2\}$ along a circle with diameter $[1, 2]$, while the remaining points $\{3, \dots,k\}$ do not move. Specifically, $\sigma_1 \in \Con^k\C$ may be represented by the loop
\[
z_1(t) = \half (3 - e^{i\pi t}),\; z_2(t) = \half(3 + e^{i\pi t}), \;\; z_j(t) = j, \; j > 2,\; 0 \leq t \leq 1.
\]
Note that $z_1(0) = z_2(1) = 1, \; z_1(1) = z_2(0) = 2$. Let us normalize and define
\[
\Gamma(z_1,\dots,z_k) := \frac{\Delta (z_1,\dots,z_k)}{|\Delta (z_1,\dots,z_k)|}
\]
Note that $\Gamma$ and $\Delta$ differ by a (variable) positive real, so we
have
\[
\Gamma_* = \Delta_* : \pi_1\Con^k\C \to \pi_1 \C^\times
\]
Basepoints in the target may differ --  for $\Gamma_*$ the basepoint is $1 \in \C^\times$ -- but this does not matter since $\pi_1 \C^\times$ is abelian. It remains to show $\Gamma_*(\sigma_1) = 1.$

Noting that if $\{z_1(t),\dots,z_k(t)\},\; 0 \leq t \leq 1$ is a loop in $\Con^k\C$ based at $\{1,\dots,k\}$ representing $\beta \in \pi_1\Con^k(\C)$, and using
\[
\Gamma(z_1,\dots,z_k) = \prod_{i < j} \frac{(z_j - z_i)^2}{|z_j - z_i|^2}
\]
then $\Gamma(z_1(t),\dots,z_k(t))$ is a loop in $\C^\times$ based at $+1$. Moreover, if we define
\[
v_{i,j}(t) := \frac{(z_j(t) - z_i(t))^2}{|z_j(t) - z_i(t)|^2} \; {\rm whenever} \; i < j
\]
we see that each $v_{i,j}(t)$ is also a loop in $\C^\times$ based at $+1$ and that $\Gamma$ is the pointwise product of these.

{\bf Lemma:}  {\it Suppose $f_m : X \to \C^\times$ are maps from a topological space $X$ to the nonzero complex numbers, and that $x_0 \in X$ satisfies $f_m(x_0) = 1$ for all $m = 1,\dots,M$. Let $f(x) = f_1(x) \cdot \dots \cdot  f_M(x)$ be the pointwise product of these complex-valued functions; then $f_*$ and $f_{m*}$ are all homomorphisms $\pi_1(X,x_0) \to \pi_1(\C^\times,1)$, and for $\alpha \in \pi_1(X,x_0)$, we have
\[
f_*(\alpha) = f_{1*}(\alpha) + \cdots + f_{M*}(\alpha) \in \pi_1(\C^\times,1) \cong \Z \;.
\]
$\Box$}

We can now calculate $\Gamma_*(\sigma_1)$ using the loop $\{z_1(t), \dots, z_k(t)\}$ described above representing $\sigma_1$.  First note that
\[
v_{1,2}(t) := \frac{(z_{2}(t) - z_1(t))^2}{|z_{2}(t) - z_1(t)|^2} =
e^{2\pi it}
\]
represents a preferred generator of $\pi_1(\C^\times, 1) \cong \Z$.  That is, $v_{1,2*}(\sigma_1) = 1$. 

It remains to show that $v_{i,k*}(\sigma_1) = 0$ if $\{i,k\} \ne \{1,2\}$, for then the lemma implies $\Gamma_*(\sigma_1) = 1$, as was to be shown. The easy case is $2<i<k\le n.$  Then $v_{i,k}(t) = 1$ is constant and therefore $v_{i,k*}(\sigma_1)=0.$ Consider the remaining cases with $i \in \{1,2\}, j > 2$.  A little trigonometry shows the unit complex number $\frac{z_j(t) - z_i(t)}{|z_j(t) - z_i(t)|}$ is restricted to have modulus in the interval
\[
[-\arcsin(1/3), \arcsin(1/3)],
\]
with $\arcsin(1/3) \sim 19.47122...$ degrees. Therefore $v_{i,j}(t)$, which is the square of this, has real part always positive. Such a loop, based at 1, is nullhomotopic in $\C^\times.$ $\qed$ 

$\; \bullet$  This corrects the junior author's false claim (in the last sentences of \cite{4}) that $\Delta_{k*} = \binom{k}{2} \w_k$. He apologizes sincerely and thanks the senior author for correcting this error, and much more. \bigskip

{\sc \S 6 Calculations in Artin coordinates}

{\bf 6.1} We have defined $\Bc_*$ to be the free left-ordered braided monoidal category on one generator. Similarly, let
 \[
\Bc_*  \supset S\Bc_* := [* \sslash S\B_*] = \coprod_{k \geq 0} [k \sslash S\B_k]
\]
be the monoidal (but no longer braided) subcategory defined by the derived subgroups of braids with writhe zero. Its morphisms can be regarded as equivalence classes, modulo type II and III (but not type I) Reidemeister moves, of isotopy classes of knot projections. 

Let $[* \sslash \Sigma_*]$ be the symmetric monoidal category of finite ordered sets and bijections, while (for a commutative ring $A$), let $[* \sslash \Gl_*(A)]$ be the category of finitely generated free $A$-modules and isomorphisms between them (sometimes abusively denoted by $(A -\Mod)$). Eventually $A$ will be the Alexander group ring $\Lambda = \Z[t,t^{-1}]$ of knot theory. 

As noted in \S 2.1, the connected component functor $\pi_0$, when applied to the morphism objects in $\Cc^*$, defines a composition 
\[
\xymatrix{
\Cc^\otimes \ar[r] &\Cc^* \ar[r]^{\pi_0} & \Bc_*}
\]
of functors which becomes a homotopy equivalence after geometric realization, and is moreover monoidal
as a functor from its first to its last term. In this coordinate system, the Burau representation $\fb : \{z\} \mapsto 
H_1(V\{z\},\Z)$ of \S 3 is a Hecke deformation\begin{footnote}{$\fb(\sigma_i)^2 = (1 - t)\fb(\sigma_i) + t \cdot {\bf 1}$}\end{footnote}
\[
\xymatrix{
[* \sslash \B_*] \ar[d]^\varepsilon \ar[r]^-\fb & [* \sslash \Gl_*(\Lambda)] \ar[d]^{t \mapsto 1} \\
[* \sslash \Sigma_*] \ar[r]^\rho & [* \sslash \Gl_* (\Z)] }
\]
of the permutation representation $\rho : S \to \Fns(S,\Z)$ of the symmetric group, with $\varepsilon$ sending a braid to its endpoint permutation. It is defined on Artin generators by 
\[
\fb(\sigma_i) :=
\left[\begin{array}{ccc}
{\bf 1}_{i-1} & 0 & 0 \\
0 & {\bf b} & 0 \\
0 & 0 & {\bf 1}_{k-1-i}
\end{array}\right] \in M_k(\Lambda) \;,
\]
where ${\bf 1}_l \in M_l(\Z)$ is the $l \times l$ identity matrix, and
\[
{\bf b} :=
\left[\begin{array}{cc}
1-t & t \\
1 & 0
\end{array}\right] \;;
\]
thus $\det \fb(\sigma_i) = -t$.\bigskip

{\bf 6.2} {\sc Group completion}

{\bf 6.2.1}  This paper is a kind of ode to group completion as the noncommutative generalization of Grothendieck's $K$-theory. 

The group completion theorem emerged from Quillen's work on algebraic $K$-theory; in the language of \cite{46} it asserts roughly that the homotopy quotient of an action by a monoid $M$ behaves in homology like an action of a group-like completion
\[
M \mapsto \Omega BM := M^+ \;.
\]
This construction sends a topological group to a homotopy equivalent $H$-space; more generally, grouplike $E_n$-spaces are precisely $n$-fold deloopings of pointed spaces, and $B$, regarded as a functor from $E_n$-spaces to grouplike $E_{n-1}$-spaces, is left adjoint to $\Omega$. In this context we have, for example, that
\[
|(A-\Mod)|^+ := |\coprod_{k \geq 0} |[k \sslash \Gl_k(A)]|^+ = \Z \times B\Gl^+_\infty(A) 
\]
defines the infinite loop-spectrum representing the algebraic $K$-theory groups $K_*(A)$ (modulo issues around $K_0(A)$). 

A useful algebraic consequence of the group completion theorem is that if $\tau \in Z(M)$ is central (and is identified with its image in the Pontryagin algebra $Z[\pi_0 M] \cong H_0(M,\Z)$), then the induced homomorphism
\[
\xymatrix{
H_*|M| \ar[r] & H_*|M|^+ }
\]
becomes an isomorphism after inverting $\tau$. For example, inverting the operation $ - \mapsto  - \otimes 1$ which adjoins an untangled strand on the right defines the homology
\[
\tau^{-1}[\bigoplus H_*(\B_k)] \cong \Z \times {\rm colim} \;  \{ H_*(\B_k) \} \cong \Z \times H_*\B_\infty
\]
of the stable braid group. 

{\bf Proposition 7} : {\it Group completion defines a morphism
\[
\xymatrix{
(|\Cc^\otimes|^+, \pp) \ar[d]^-{\bpi \; \simeq} \ar[r]^-{|\fb|^+} & |(\Lambda-\Mod)|^+, \oplus) \ar[d]^\simeq \\
(|\Bc|_*,\otimes) \ar@{.>}[r]  & (\Omega^\infty K(\Lambda),+) }
\]
of monoids, factoring the homotopy-theoretic Burau representation through its classical algebraic image.} 

The Burau construction thus defines an interesting class in $K(\Lambda)^0(\Omega^2 S^3)$ \cite{61}. A topological quantum field theory, in the sense of Segal, Atiyah, and others, is (very roughly) a monoidal functor from a category of geometric objects, with something like manifolds as objects, cobordisms as morphisms, and sometimes (\eg topological gravity) with diffeomorphisms as 2-morphisms, to a more algebraic (often symmetric monoidal) category. In this broad sense the Burau representation $\Cc^\otimes \to [* \sslash \Gl_*(\Lambda)]$ (or Whitehead torsion \cite{47} for that matter) is a TQFT; but the Burau image  $\fb_*(c^\Bc_{k,k'})$ of the Joyal-Street twist is {\bf not} the Koszul twist of homological algebra, so the Burau representation does {\bf not} define a {\bf braided} homotopical field theory.

{\bf 6.2.2} The harmonic\begin{footnote}{\ie $\Delta \log p\{z\} = 4\pi \sum \delta(t - z_i)$: the Hodge/Laplace operator $\Delta$ should not be confused with the discriminant.}\end{footnote} potential $\log p\{z\}(t)$ of Newton and Segal, where
\[
t,\{z\} \in \C \times \Con^*(\C) \mapsto p\{z\}(t) := \prod (t - z_i) \in \C[t] \;,
\] 
defines a one-form
\[
(d \log p)\{z\} := p(t)^{-1} dp(t)  = \sum (t - z_i)^{-1} dt \in \Omega^1\Con^*(\C) \;,
\] 
and regarding $\C$ as a Riemann surface parametrized by $t$ defines a family
\[
\{z\} \mapsto (p^{-1} \nabla p)(t) : (\C_+,+) \to (\C_+,0)
\]
of gradient vector fields (beware basepoint change), a kind of statistically mechanical physical field
\[
\xymatrix{
\nabla \log p_* :  \Con^*(\C) \ar[r]^-\simeq & \Omega^2_{* \geq 0} S^2 }
\]
defined away from the discriminant locus $\Delta^{-1}(0) \subset \C^k/\Sigma_k$. Following \cite{55}(\S3), it is homotopic to the group completion map. 

The integral cohomology of the pure braid space $\OCon^k(\C)$ is torsion-free\begin{footnote}{When $n>2$ the homology and cohomology of $\OCon^k(\R^n)$ are algebraically similar \cite{27,28}, but with generators in degree $n-1$}\end{footnote}, with $H^*_\dR(\Con^k(\C))$ of rank $\binom{n}{2}$ generated by the classes of Arnol'd's closed one-forms
\[
\alpha_{ij} = \alpha_{ji} = \frac{dz_i - dz_j}{z_i - z_j} \in \Omega^1(\OCon_k(\C))
\]
($1 \leq i \neq j \leq k$), subject to the quadratic relations
\[
\alpha_{ij}\alpha_{jk} +\alpha_{jk}\alpha_{ki} + \alpha_{ki}\alpha_{ij} = 0 \;.
\]
The symmetric group acts on the de Rham cohomology
\[
H^0(\Sigma_k;H^*_\dR(\OCon^k(\C))) \cong H^*_\dR(\Con^k(\C)) 
\]
identifying $\Con^k(\C)$ as a rational homology circle. The discriminant (\S 2.3)
\[
\Delta^*_k : H^*_\dR(\C^\times) \to H^*_\dR(\Con^k(\C))
\]
defines an isomorphism sending the Cauchy generator $\fc := (2\pi iz)^{-1}dz \in \Omega^1(\C^\times)$ to
\[
\Theta^*(\fc) = \Delta_k^*(\fc) = (2 \pi i \Delta_k)^{-1} d \Delta_k  = (2 \pi i)^{-1} \sum_{i \neq j} \alpha_{ij} \;,
\]
presenting the discriminant as a generalized Abel-Jacobi map. We are indebted to N Kitchloo for pointing out that $\Theta$ is a connection form for the Knizhnik-Zamolodchikov line bundle over the configuration space \cite{38} associated to a representation of the Lie algebra $\gl$ of the multiplicative group (as for example in \cite{12}).
 
{\bf Proposition 8} : {\it The elements $d \log p_*$ and $\Theta_*(\fc)$ of $\Omega^1(\Con^k(\C))$ are cohomologous.}

{\bf Proof:} Composing the adjoint 
\[
S^1 \cong \R_+ \to \Maps_*(\R^2_+,\R^3_+) \cong \Omega^2 S^3
\]
to stereographic projection ($\xymatrix{\R_+ \wedge \R^2_+ \ar[r]^-\cong & \R^3_+}$) with a degree $k$ self-map of the circle $\R_+ \cong \T$ defines a generator $h_k \in \pi_1(\Omega^2_k S^3)$ with Hurewicz image Kronecker dual to 
a generator of $H^1(\Omega^2_k S^3) \cong \Z$. The two maps $\nabla \log_* p$ and $h_* \Delta_*$ both induce the identity isomorphism on fundamental groups and thus on $H^1_\dR$, so their difference on a generator must be exact. \bigskip

{\bf 6.3} The theorem of Gorin and Lin defines a fibration
\[
\xymatrix{
{\tC^* = \coprod_{k \geq 0} [\tCon^k(\C) \sslash \D]} \ar[r] & \Cc^* \ar[r]^-\Delta & [\N \sslash \T]}
\]
of topological categories, resulting in an exact sequence
\[ 
\bpi |\tC|^* = {\rm S}\Bc_* \to \bpi |\C|^* = \Bc_* \to [\N \sslash \Z]
\]
of fundamental groupoids, and similarly for $|\tC|^\otimes$. The covering space is homotopy equivalent to the Milnor fiber of the discriminant. 
 
 {\bf Proposition 9} : {\it Group completion defines a homotopy-commutative diagram
\[
\xymatrix{
\Z \ar[r] & |\tC_\otimes|^+ \ar[d]^-\simeq \ar[r] & |\Cc_\otimes|^+ \ar[d]^-\simeq \\
{} & \Z \times \Omega^2 S^3 \langle 3 \rangle \ar[r] & \Z \times \Omega^2 S^3 \ar[r] & \Z \times \T }
\]
of topological monoids.} 

Proof: The last line in the diagram splices the equivalence
\[
\Z \times \Omega^2 S^3 \simeq \Omega^2_* S^2
\]
(defined by taking the double loops on the Hopf map, resulting in a semi-direct product of group-like spaces), with double loops on the 3-connected cover 
\[
H(\Z,2) \to S^3 \langle 3 \rangle \to S^3 
\]
of the three-sphere. The right-hand vertical equivalence goes back to Boardman and Segal, while that on the left is the subject of the appendix to \cite{4}. $\Box$ \bigskip

{\bf Remarks} Following \cite{59}(Th 1, 3b), delooping $|\Cc_*|^+ := \Omega B (|\Cc|_*) \simeq \Omega^2_*S^2 $ implies an equivalence $B (|\Cc|_*) \simeq \Omega S^2$, with a universal cover $B (|\tC|_*) \simeq \Omega S^3$ which inherits an action of $\Z$ by deck-translations. The Pontryagin algebra $H_*\Omega S^2$ is commutative polynomial on a generator of degree {\bf one}, while $\Omega S^3 \cong (\Omega \Sigma) (S^2)$ is stably equivalent to James' free topological monoid on $S^2$, with polynomial homology on a generator of degree two. The loop group $\Omega S^3 \cong \Omega \SU(2)$ has a remarkable level one projective representation which embeds its homology algebra 
\[
H_*\Omega \SU(2) \cong \Z[b] \to \Z[\gamma_k(b) \:|\: k \geq 1] \cong H_*\C P^\infty
\]
in its divided power envelope by  $\gamma_k(b) \mapsto b^k/k!$.  The deck-transformation group acts by changing the level structure of the representation.

The space $\Omega^2 S^3$ of smooth maps from $(\R^2_+,+)$ to $(\SU(2),1)$ is the state space $\wz$ for the simplest Wess-Zumino-Witten model in string theory, where the infinite-dimensional 3-connected groups Spin$(n)\langle 3 \rangle$ are sometimes \cite{48} called $\String(n)$. This proposition might thus be paraphrased as asserting that the universal covering space $\wzw := \Omega^2 \String(3)$ is a homotopy torsor for a monoid action of $|\tC|^\otimes$, defining a particle model for 2D fields on $\R^3_+$ with cubic (Chern-Simons-Cartan-Witten) interactions. \bigskip

{\bf \S 7} {\sc Monoidal braidings and $p$-adic spheres} \bigskip

{\bf 7.1} A commutativity constraint $c^\A$ (\eg in the sense of \cite{37}(Ch XIII)) on a monoidal category $\A$ is a natural equivalence of the product functor $\otimes : A \times \A \to \A$ with the twisted product $\otimes \circ \tau$:
\[
\xymatrix{
\A \times \A \ar[r]^\tau & \A \times \A \ar[r]^\otimes & \A \;, }
\]
where $\tau$ is the flip functor, defined on objects by $\tau(X,Y) = (Y,X)$. The constraint is required to respect identity elements and to satisfy certain hexagon axioms, which make $\A$ braided monoidal.

The category $\Bc_*$ is braided by compositions with elements 
 \[
 c^\Bc_{k,k'} := \sigma_{[k:1]} \cdots \sigma_{[k+k'-1:k']} 
 \]
(where $\sigma_{[k:l]} = \sigma_k \cdots \sigma_l, \; k \geq l$) which pass a trivial $k$-strand braid on the left over a trivial $k'$-strand braid on the right, \eg $c^\B_{2,2} = \sigma_2 \sigma_1 \sigma_3 \sigma_2$. The image 
\[
\varepsilon \circ c^\Bc_{k,k'} := c^\Gl_{k,k'}  = 
\left[\begin{array}{cc}
0 & {\bf 1}_k \\
{\bf 1}_{k'} & 0
\end{array}\right]  \in \Sigma_{k+k'} \subset \Gl_{k+k'}(\Z)
\]
of this Joyal-Street constraint under endpoint permutation is the usual Koszul braiding of graded algebra.

The writhe of a braid thus equals the determinant of its Burau representative: 

{\bf Proposition 10} : {\it The functor
\[
\det \fb : [* \sslash \B_*] \to [* \sslash \Lambda^\times]
\]
sends a braid $b$ to} $(-t)^{\w(b)}$, \eg \;  $(\det \fb)(c^\Bc_{k,k'}) = (-t)^{kk'}$. 

 {\bf Definition} Given a commutative ring $A$, define the Picard category $\Pic_@(A)$  of anyonic $A$-lines \cite{44} to be the groupoid with free rank one $A$-modules $A[k], \; k \in \Z$ as objects, and morphisms
\[
\mor(A[k],A[k']) = A^\times \; {\rm if} \; k = k', \; = \emptyset \; {\rm otherwise} \;,
\]
together with the monoidal structure 
\[
A[k],A[k'] \mapsto A[k] \otimes_A A[k'] = A[k+k'] \;.
\]
The choice of an element $a \in A^\times$ defines a commutativity constraint 
\[
\xymatrix{
c^\A_{k,k'} : A[k] \otimes_A A[k'] \ar[r] & A[k+k'] \ar[r]^-{a^{kk'} \cdot } & A[k'+k] \ar[r] & A[k'] \otimes_A A[k]}
\]
(which makes the diagram 
\[
\xymatrix{
X \otimes_A Y \ar[d]^{\alpha \otimes \beta} \ar[r]^{c^\A_{X,Y}} & Y \otimes_A X \ar[d]^{\beta  \otimes \alpha} \\ 
X' \otimes_A Y' \ar[r]^{c^\A_{X',Y'}} & Y' \otimes_A X'  \;.} 
\]
commutative, for any $\alpha, \; \beta \in \mor(X,X'), \; \mor(Y,Y')$). 

The category $\Pic_@(A)$ is equivalent to the groupoid with graded free rank one $A$-modules as objects, and isomorphisms of such things as morphisms. Its geometric realization ${\rm Pic}_@(A)$ is a two-stage Postnikov system with $\pi_0 = \Z$ and $\pi_1 = A^\times$ with Postnikov invariant 
\[
[\beta : k \mapsto a^k] \in H^1(\Z,A^\times) \;.
\] 
In this language the writhe of a braid becomes the top exterior power of Burau as a $\Lambda$-module. The group of equivalence classes of such $\Lambda$-lines is 
\[
\pi_0|\Pic(\Lambda)| \cong K_1(\Lambda) \cong \Lambda^\times = \{\pm t^w \:|\; \:w \in \Z \} \;,
\]
for example by a theorem of Higman on units of group rings, or as an instance $K_*(\Z[t,t^{-1}]) \cong K_*(\Z) \oplus K_{*-1}(\Z)$ of the fundamental  theorem of algebraic $K$-theory \cite{60}(Ch V).

 {\bf Proposition 11} : 
 \[
 \w : b \mapsto \det \fb (b)  : [* \sslash \B_*] \to \Pic_@(\Lambda)
\]
{\it is a {\bf braided} \cite{16} monoidal functor}. 

The composition $\xymatrix{\Cc^\otimes \ar[r]^-{\pi_0} & \Bc_* \ar[r]^-\w & \Pic_@(\Lambda)}$ can be interpreted as an invertible anyonic TQFT for finite subsets of the plane, with diffeomorphisms of their complements as morphisms, see below. It is as braided as it can be, given that its domain is braided only modulo higher homotopies. 

{\bf 7.2} In the language of stable homotopy theory, a $k$-dimensional $p$-adic sphere $X$ is a $p$-complete spectrum with reduced homology $\tilde{H}_*(X,\Z)$ free of rank one over $\Z_p$ in degree $k$ \cite{11},\cite{59}(\S 3); it is oriented by a choice $\Z_p[k] \cong \tilde{H}_*(X,\Z)$ of generator. (Oriented) $p$-adic spheres are invertible under (completed) smash product in the category of $E_\infty$ modules over the stable $p$-adic sphere ring-spectrum $S^0_p$. The Picard category of such things, with homotopy equivalences as morphisms, has  geometric realization
\[
|\Pic^\Or(S^0_p)| \simeq \Z \times B\Sl_1(S^0_p)
\]
homotopy equivalent to the identity component of the $p$-completion of the classifying space of the monoid $\SF(\infty)$ of degree one self-maps of the stable sphere \cite{41}. An element $a \in \Z^\times_p$ defines a commutativity constraint and a braided monoidal category $\Pic^\Or_@(S^0_p)$. 

Recent work \cite{5}(Prop 5.9), \cite{18} of broad interest defines analogs of Haar measure for $p$-adic Lie groups in terms of $p$-adic spheres 
\[
\ess(\g) := {\rm hocolim}_{i \to \infty} \Sigma^\infty_+(Bp^i\g), \; \tilde{H}_* (\ess(\g),\Z) \cong \Lambda^{\rm top} \g
\]
associated to their Lie algebras $\g$, the maps of the system being stable traces (defined using Dyer-Lashof algebras as in \cite{50}(\S 3) (which Pontryagin duality simplifies in this abelian case)).

{\bf Proposition 12} : {\it The (compact, commutative) $p$-adic Lie algebras $H_1(V\{z\},\Z_p)$ define a commutative diagram
\[
\xymatrix{
\Cc^\otimes \ar[d]^\bpi \ar@{.>}[r]^-{\ess \fb} & \Pic^\Or_@(S^0_p) \ar[d] \\
\Bc_* \ar[r] & \Pic_@(\Z_p) }
\]
lifting the Burau representation to a braided monoidal functor with values in the category of oriented $p$-adic spheres.} 

A ring homomorphism $t \mapsto a \in \Z_p^\times : \Lambda \to \Z_p$ defines a monoidal specialization functor $\Pic_@(\Lambda) \to \Pic_@(\Z_p)$, sending $t$ to an element of the topologically cyclic (\cf Atiyah-Tall or Serre) extension
\[
\xymatrix{
1 \ar[r] & \mu_{p-1} \ar[r] & \Z_p^\times \ar[r] & (1 + p\Z_p )^\times \ar[r]^-{\log \; \cong} & \Z_p \ar[r] &1 } 
\] 
($p$ odd). Any element $1 + kp$ with $k \neq 0$ mod $p$, eg $1 - p$, generates $\Z_p^\times$. Up to regrading and rescaling, then, the braiding on $\Pic^\Or_@(S^0_p)$ is essentially unique. 

Group completion now defines a commutative diagram
\[
\xymatrix{
|\Cc^\otimes|^+ \ar[d]^-\simeq \ar[r]^-{|\ess \fb|^+} & {\rm Pic}(S^0_p) \ar[d] \\
|\Bc_*|^+ \ar[r] \simeq \Z \times \Omega^2 S^3 \ar[r] & {\rm Pic}(\Z_p) }
\]
of group-like spaces. In dimensions  0 and 1, Clausen's $p$-adic $J$-homomorphism \cite{18}(\S 3.2 4.3) 
\[
\xymatrix{
K(\Z_p) \ar[r]^-{J_{\Z_p}} & {\rm Pic}(S^0_p)  \ar[r] & {\rm Pic}(L_{K(1)}S^0) \ar[r]^-{R\log_{K(1)}} & \Sigma L_{K(1)}S^0 \simeq 
\Sigma \; {\rm image} \; J_p \;,}
\]
composed with Rezk's logarithmic operation \cite{52}, is an isomorphism, implying that $|\ess \fb|^+$ lifts to a morphism 
\[
|\tC^\otimes|^+ \simeq |S\Bc_*|^+ \simeq \Omega^2 S^3 \langle 3 \rangle  \to \Sigma \; {\rm image}\; J
\]
of spectra. The constructions of \cite{2}(\S 2.5) define $p$-adic Thom spectra as colimits of such families of $S^0_p$ lines; patching them together defines a global object which we will call $M\SF(3)$. \bigskip

{\bf \S 8} {\sc generalized Steenrod representability}

F. Cohen \cite{20} showed, by considering the composition $\B_k \to \Sigma_k \to {\rm O}(k)$, that every mod two homology class can be represented by a map from a manifold with a lift of its tangent structure from the orthogonal group to a braid group. The situation at odd primes is considerably deeper, involving \cite{21}(Prop 9) group completion. The theorem of Mahowald and Hopkins \cite{42, 43} implies, following Quillen's conventions, that a proper map $V \to X$ between Poincar\'e duality spaces \cite{13,40,51}, oriented by a lift
\[
\nu_{V/Z} : V \to \Omega^2 S^3\langle 3 \rangle \simeq B\SF(3) \simeq \wzw
\]
of its spherical normal fibration from $B\SF$, represents an integral homology class in $X$, and conversely. Such matters are of interest in geometric nonlinear analysis \cite{28}(\S 3.7,5.11,9.6), \cite{30}. 

In 1967 Peterson and Toda \cite{49} showed that the spectrum $M\SF$ (defined \cite{41}(Ch III) by spherical Thom fibrations) is generalized Eilenberg - Mac Lane. One consequence of the 

{\bf Theorem of Hopkins and Mahowald}
\[
H\Z \to M\SF(3) \simeq M\wzw \to M\SF
\]
{\it makes $M\SF$ an $E_2$-module ring spectrum over $H\Z$} 

is that, away from the prime two, Poincar\'e duality cobordism is decided by spherical characteristic numbers, and thus by classical (\eg torsion) characteristic numbers. \bigskip

The proof in \cite{21}(lemma 7), \cite{9}(\S 9) (see also \cite{1, 45}) characterizes maps $\Omega^2 S^3 \langle 3 \rangle  \to B\SF$ which pull back a Thom spectrum homotopy equivalent to $H\Z$ in terms of an invariant with values in $\pi_{2(p-1)}B\SF \cong \pi^S_{2p-3}$. The Burau determinant is a double loop map; this structure is respected by the action of the first Dyer-Lashof operation. Applying that operation to the class in dimension one defined by the discriminant defines a class which generates the first (Adams) class $\alpha_1$ in the image of the $J$-homomorphism. $\Box$
 
In \cite{19}(Ch IV) Cohen moreover proved that the Pontryagin rings $H_*\SF(n)$ are commutative, which suggests that an extension
\[
\SF(3) \to \SF \to \SF / \SF(3)  \simeq \SF_{h\SF(3)}
\]
of Rognes-Hopf-Galois \cite{53} objects with
\[
M\SF \wedge_{H\Z} M\SF \simeq S^0[\SF/\SF(3)] \wedge_{H\Z}  M\SF \;;
\]
and an accessible Rothenberg-Steenrod spectral sequence
\[
 \Tor^{H_*\SF(3)}_* (H_*\SF,\Z)  \Rightarrow H_*(\SF / \SF(3)) 
\]
makes sense. \bigskip

\bibliographystyle{amsplain}

\end{document}